\documentclass[a4paper,11pt]{article}
\usepackage{graphicx}
\usepackage{amsfonts}
\usepackage{amsmath}
\usepackage{amssymb}
\usepackage{indentfirst}
\usepackage{titlesec}
\usepackage{caption2}
\usepackage{color}

\def\barr{\begin{array}}
\def\earr{\end{array}}
\def\bali{\begin{aligned}}
\def\eali{\end{aligned}}
\def\bearr{\begin{eqnarray}}
\def\eearr{\end{eqnarray}}

\providecommand{\play}{\displaystyle}

\providecommand{\li}{\limits}

\providecommand{\pt}{\partial}

\providecommand{\ra}{\rightarrow}

\providecommand{\da}{\downarrow}

\providecommand{\Prob}{\mathbf P}

\providecommand{\E}{\mathbf E}

\providecommand{\al}{\alpha}

\providecommand{\bt}{\beta}

\providecommand{\gm}{\gamma}

\providecommand{\Gm}{\Gamma}

\providecommand{\dt}{\delta}

\providecommand{\ve}{\varepsilon}

\providecommand{\tht}{\theta}

\providecommand{\kp}{\kappa}

\providecommand{\lb}{\lambda}

\providecommand{\sm}{\sigma}

\providecommand{\N}{\mathbb N}

\providecommand{\R}{\mathbb R}

\providecommand{\cF}{\mathcal F}

\providecommand{\iY}{\mathfrak{Y}}

\providecommand{\iiY}{\mathbf{\mathfrak{Y}}}

\providecommand{\1}{\mathbf 1}

\providecommand{\contfunc}{\mathbf{C}}

\providecommand{\boldx}{\boldsymbol{x}}

\providecommand{\boldX}{\boldsymbol{X}}

\providecommand{\boldn}{\boldsymbol{n}}

\providecommand{\vphi}{\varphi}

\providecommand{\ubar}{\bar{u}}

\providecommand{\Vbar}{\overline{V}}

\begin{document}

\title{Wave front propagation for a reaction-diffusion equation in narrow random channels}
\author{Mark Freidlin\thanks{Department of Mathematics,
University of Maryland at College Park, mif@math.umd.edu.} \ , \
Wenqing Hu\thanks{Department of Mathematics, University of Maryland
at College Park, huwenqing@math.umd.edu.}}

\date{
}

\maketitle

\begin{abstract}

We consider a reaction-diffusion equation in narrow random channels.
We approximate the generalized solution to this equation by the
corresponding one on a random graph. By making use of large
deviation analysis we study the asymptotic wave front propagation.

\end{abstract}

\textit{Keywords}: reaction-diffusion equation, wave front
propagation, diffusion processes on graphs, random environment.

\textit{2010 Mathematics Subject Classification Numbers}: 35K57,
35A18, 60J60, 60K37.

\section{Introduction}

In studying the motion of molecular motors we introduced in
\cite{[F-H molecular motors]} a solvable model: we think of the
molecular motors as diffusion particles traveling in a narrow random
channel. Based on the model suggested in \cite{[F-H molecular
motors]}, we consider in this paper wave front propagation for a
reaction-diffusion in narrow random channels. Problems of this type
naturally appear in the theory of nerve impulse propagation and in
combustion theory. Our analysis relies on techniques in large
deviations similar to that of \cite[Chapter 7]{[F red book]} and
\cite{[Nolen-Xin KPP random drift]}, \cite{[Nolen-Xin CMP]},
\cite{[Taleb]}, \cite{[Zeitouni LDP RWRE]}, \cite[Chapter 5]{[Xin
book]}. We shall note that problems of this type are mentioned in
\cite[Chapter 7]{[F green book]}, \cite{[F-Sheu]},
\cite{[F-Kostas]}. It is also interesting to note that similar
problems are considered in \cite{[Molchanov-Vainberg]},
\cite{[Lions-Souganidis CPDE]}, \cite{[Lee et al]} but from
different points of view.

\begin{figure}
\centering
\includegraphics[height=7cm, width=10cm , bb=24 20 403 289]{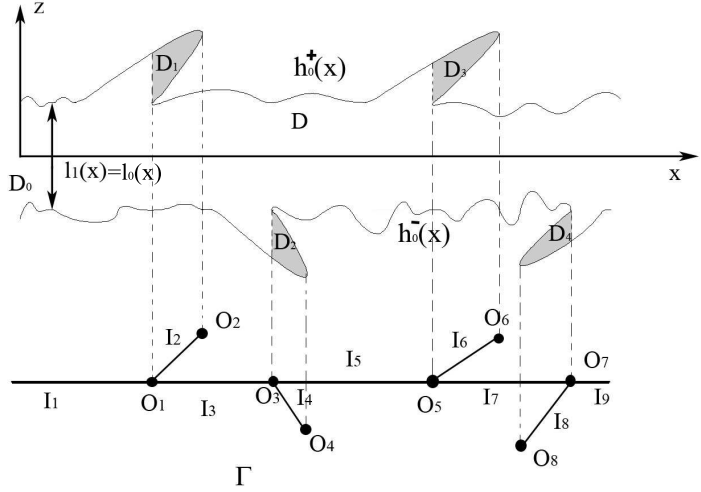}
\caption{A model of the molecular motor.}
\end{figure}

Let us first briefly recall the model introduced in \cite{[F-H
molecular motors]}. Let $h_0^{\pm}(x)$ be a pair of piecewise smooth
functions with $h_0^+(x)-h_0^-(x)=l_0(x)>0$. Let $D_0=\{(x,z): x\in
\R , h_0^-(x)\leq z \leq h_0^+(x)\}$ be a tubular $2$-d domain of
infinite length, i.e. it goes along the whole $x$-axis. At the
discontinuities of $h_0^{\pm}(x)$, we connect the pieces of the
boundary via straight vertical lines. The domain $D_0$ models the
"main" channel in which the motor is traveling. Let a sequence of
"wings" $D_j$ ($j\geq 1$) be attached to $D_0$. These wings are
attached to $D_0$ at the discontinuities of the functions
$h_0^{\pm}(x)$.

Consider the union $D=D_0\bigcup\left(\bigcup\li_{j=1}^\infty
D_j\right)$. An example of such a domain $D$ is shown in Fig.1, in
which one can see four "wings" $D_1,D_2,D_3,D_4$. We assume that,
after adding the "wings", for the domain $D$, the boundary $\pt D$
has two smooth pieces: the upper boundary and the lower boundary.
Let $\boldn(x,z)=(n_1(x,z), n_2(x,z))$ be the inward unit normal
vector to $\pt D$. We make same assumptions as in \cite{[F-H
molecular motors]}.

\textbf{Assumption 1.} The set of points $x\in \R$ for which there
are points $(x,z)\in \pt D$ at which the unit normal vector
$\boldn(x,z)$ is parallel to the $x$-axis: $n_2(x,z)=0$ has no limit
points in $\R$. Each such point $x$ corresponds to only one point
$(x,z)\in \pt D$ for which $n_2(x,z)=0$.

\textbf{Assumption 2.} For every $x$ the cross-section of the region
$D$ at level $x$, i.e., the set of all points belonging to $D$ with
the first coordinate equal to $x$, consists of either one or two
intervals that are its connected components. That is to say, in the
case of one interval this interval corresponds to the "main channel"
$D_0$; and in the case of two intervals one of them corresponds to
the "main channel" $D_0$ and the other one corresponds to the wing.
The wing will not have additional branching structure. Also, for
some $0<l_0<\bar{l}_0<\infty$ we have $\l_0\leq
h_0^+(x)-h_0^-(x)=l_0(x)\leq \bar{l}_0$.

Let us take into account randomness of the domain $D$. Keeping the
above assumptions in mind, we can assume that the functions
$h_0^{\pm}(x)$ and the shape of the wings $D_k$ ($k=1,2,...$) are
all random. Thus we can view the shape of $D$ as random. We
introduce a filtration $\cF_s^t$, $-\infty\leq s<t\leq \infty$ as
the smallest $\sm$-algebra corresponding to the shape of
$D\cap\{(x,z): x\in [s, t]\}$. We introduce stationarity and mixing
assumptions. Let us consider some $A\in \cF_s^t$, $-\infty \leq
s<t\leq \infty$. The set $A$ consists of some shapes of the domain
$D\cap \{(x,z): x\in [s,t]\}$. Let $\tht_r$ ($r\in \R$) be the
operator corresponding to the shift along $x$-direction:
$\tht_r(A)\in \cF_{s+r}^{t+r}$ consists of the same shapes as those
in $A$ but correspond to the domain $D\cap \{(x,z): x\in
[s+r,t+r]\}$.

\textbf{Assumption 3.} (stationarity) We have
$\Prob(A)=\Prob(\tht_r(A))$.

\textbf{Assumption 4.} (mixing) For any $A\in \cF_s^t$ and any $B\in
\cF_{s+r}^{t+r}$ we have
$$\lim\li_{r\ra \pm\infty}\sup\li_{A\in \cF_s^t, B\in \cF_{s+r}^{t+r}}|\Prob(A\cap B)-\Prob(A)\Prob(B)|=0$$
exponentially fast.

For instance, we can assume that there exists some $M>0$ such that
$\Prob(A\cap B)=\Prob(A)\Prob(B)$ for $|r|\geq M$.

In particular, the mixing assumption implies that the transformation
$\tht_r$ is ergodic.

Here and below the symbols $\Prob$ and $\E$ etc. refer to
probabilities and expectations etc. with respect to the filtration
$\{\cF_s^t\}_{-\infty\leq s <t\leq \infty}$.

Let $D^\ve=\{(x,\ve z): (x,z)\in D\}$. The parameter $\ve>0$ is
small. The domain $D^\ve$ models the narrow random channel. Let us
consider the following reaction-diffusion equation in the domain
$D^\ve$:

$$\left\{\begin{array}{l}
\dfrac{\pt u^\ve}{\pt t}=\dfrac{1}{2}\left(\dfrac{\pt^2 u^\ve}{\pt
x^2} +\dfrac{\pt^2 u^\ve}{\pt z^2}\right)+V(x,z)\dfrac{\pt
u^\ve}{\pt x}+ f(u^\ve) \ ,
\\
u^\ve(0,x,z)=g(x) \ ,
\\
\left.\dfrac{\pt u^\ve}{\pt \nu}\right|_{\pt D^\ve}=0 \ ,
\\
u^\ve=u^\ve(t,x,z) \ , \ (t,x,z)\in \R_+\times D^\ve \ . \end{array}
\right. \eqno(1.1)$$

Here $\nu$ is the inward unit normal vector field on $\pt D^\ve$;
$V(x,z)$ is the velocity field; the function $f(u^\ve)$ is smooth
and is of KPP type: $f\in \contfunc^{(\infty)}([0,1])$, $f(0)=f(1)$,
$0\leq f(u^\ve)\leq f'(0) u^\ve$ for all $u^\ve\in (0,1)$, e.g.
$f(u^\ve)=u^\ve(1-u^\ve)$. The initial function $g(x)\geq 0$ (not
identically equal to $0$) is smooth and compactly supported in $x\in
\R$: $g(x)\in \contfunc^{(\infty)}_0(\R)$. We notice that the
initial function $g(x)$ depends only on the variable $x\in \R$ and
is independent of $z$.

Alternatively, problem (1.1) can be considered on the domain $D$
with a change of variable $z\mapsto z/\ve$. The equivalent problem
takes the form

$$\left\{ \begin{array}{l}
\dfrac{\pt u^\ve}{\pt t}=\dfrac{1}{2}\left(\dfrac{\pt^2 u^\ve}{\pt
x^2} +\dfrac{1}{\ve^2}\dfrac{\pt^2 u^\ve}{\pt
z^2}\right)+V(x,z)\dfrac{\pt u^\ve}{\pt x}+ f(u^\ve) \ ,
\\
u^\ve(0,x,z)=g(x) \ ,
\\
\left.\dfrac{\pt u^\ve}{\pt \nu^\ve}\right|_{\pt D}=0 \ ,
\\
u^\ve=u^\ve(t,x,z) \ , \ (t,x,z)\in \R_+\times D \ . \end{array}
\right. \eqno(1.2)$$

Here $\nu^\ve=(\nu^\ve_1, \nu^\ve_2)=(\ve n_1(x,z), n_2(x,z))$ is
the inward unit co-normal vector field on $\pt D$ corresponding to
the operator $\dfrac{1}{2}\left(\dfrac{\pt^2 }{\pt x^2}
+\dfrac{1}{\ve^2}\dfrac{\pt^2}{\pt z^2}\right)$.

The diffusion process $\boldX_t^\ve=(X_t^\ve, Z_t^\ve)$
corresponding to problem (1.2) takes the form of equation (2) in
\cite{[F-H molecular motors]}. We have

$$\left\{\begin{array}{l}
dX_t^\ve=dW_t^1+V(X_t^\ve, Z_t^\ve)dt+\nu_1^\ve(X_t^\ve,
Z_t^\ve)d\ell_t^\ve \ ,
\\
dZ_t^\ve=\dfrac{1}{\ve}dW_t^2+\nu_2^\ve(X_t^\ve, Z_t^\ve)d\ell_t^\ve
\ ,
\end{array} \right. $$

Let $\Prob^W$, $\E^W$ denote probabilities and expectations with
respect to the filtration generated by $\boldX_t^\ve$ (and
henceforth $(W_t^1, W_t^2)$). We note that as in \cite{[F-H
molecular motors]} the motion of $\boldX_t^\ve$ is independent of
the random shape of $D$ (and henceforth $D^\ve$). We have, in the
same way as \cite{[F-H molecular motors]}, the following.

\textbf{Assumption 5.} The process $(W_t^1, W_t^2)$ is independent
of the filtration $\{\cF_s^t\}_{-\infty\leq s<t\leq \infty}$
corresponding to the shape of $D$.

We shall also make some assumptions parallel to Assumptions 7 and 8
in \cite{[F-H molecular motors]}. To this end we let $L$ be the
random variable distributed the same as the distance along $x$-axis
between two wings: $L$ is the distance along $x$-axis between two
cross-sections of $D$ where there is a branching. Let
$l_{\text{wing}}$ be the cross-section width of the wing. Let $r$ be
the projection length of a wing onto $x$-axis ($r$ can be positive
or negative; compare with \cite[Assumptions 7 and 8]{[F-H molecular
motors]}). We assume the following.

\textbf{Assumption 6.} With $\Prob$ probability $1$ we have (1)
$+\infty>\overline{L}\geq L\geq \underline{L}>0$ for some constants
$\overline{L}, \underline{L}>0$; (2) $0\leq l_{\text{wing}}(x)\leq
A_1$ and $|r|\leq A_1$ for some constant $A_1>0$.

\

Our goal in this paper is to study the asymptotic wave front
propagation properties for the generalized solution of (1.2). To be
precise, by a \textit{generalized solution} of (1.2) we mean the one
defined via the path integral representation (Feynmann-Kac) formula:

$$u^\ve(t,(x,z))=\E_{(x,z)}^W\left[\exp\left(\int_0^t
c(u^\ve(t-s, \boldX_s^\ve))ds\right)g(X_t^\ve)\right]  \ .
\eqno(1.3)$$

Here $c(u)=\dfrac{f(u)}{u}$ for $u>0$ and $c(0)=\lim\li_{u\da
0}\dfrac{f(u)}{u}=\sup\li_{u>0}\dfrac{f(u)}{u}$. The latter equality
is due to the KPP nonlinearity assumption. We shall also suppose
that $|c'(u)|\leq \text{Lip}(c)<\infty$, $u\in [0,1]$. The proof of
existence, uniqueness and regularity of the generalized solution to
the integral equation (1.3) is close to \cite[Chapter 5, Section
3]{[F red book]}. For the reader's convenience we will prove it in
Section 3 of this paper.

We introduced in \cite{[F-H molecular motors]} the metric graph
$\Gm$ corresponding to the domain $D$ (see Fig.1). Let $Y_t=(X_t,
k_t)$ be defined on $\Gm$ as in \cite[Section 2]{[F-H molecular
motors]}. The construction of the graph $\Gm$ and the process $Y_t$
as well as some basic convergence results in \cite{[F-H molecular
motors]} will be recalled in Section 2. The process $Y_t$ is a
diffusion process on $\Gm$ with a generator $A$ and the domain of
definition $D(A)$. We consider the reaction-diffusion equation
associated with the Markov process $Y_t$. This equation takes the
form

$$\dfrac{\pt u}{\pt t}=A u + f(u) \ , \ u(0,(x,k))=g(x) \ ,
\ u=u(t,y), \ (t,y)\in \R_+\times \Gm \ . \eqno(1.4)$$

The initial function $g=g(x)$ is the same as in (1.2). We require
that for each fixed $t\geq 0$, $u(t, \bullet)\in D(A)$. This
requirement ia a kind of boundary condition. For details we refer to
\cite{[F-H molecular motors]} and \cite[Chapter 8]{[FW book]} and
the references therein.

The generalized solution to (1.4) is defined as the solution to the
integral equation

$$u(t,(x,k))=\E_{(x,k)}^W\left[\exp\left(\int_0^t c(u(t-s, Y_s))ds\right)g(X_t)\right] \ . \eqno(1.5)$$

The existence and regularity of the solution can be proved in a same
way as those for (1.3). We will briefly mention this in Section 3.

We will show, in Section 3 of this paper, that as $\ve\da 0$, the
solution $u^\ve(t,(x,z))$ of (1.3) will converge to $u(t,(x,k))$ of
(1.5) in the strong sense. Here $(x,k)=\iiY((x,z))$ and the mapping
$\iiY:D\ra \Gm$ is an identification map that will be recalled in
Section 2.

After we get convergence results we will focus on the study of the
solution $u(t,(x,k))$ of (1.5). We will show that, as $t\ra\infty$,
the solution $u(t,(x,k))$ behaves asymptotically as a traveling
wave. This wave is traveling in both positive and negative
directions along $x$-axis.

The paper is organized as follows: in Section 2 we recall some
necessary basic set up of \cite{[F-H molecular motors]}, such as the
construction of the graph $\Gm$, the process $Y_t$, etc.; in Section
3 we show the convergence as $\ve\da 0$ of the solution
$u^\ve(t,(x,z))$ to $u(t,(x,k))$; in Section 4 we obtain some
auxiliary results that will be used in later sections; in Section 5
we derive the large deviation principle; in Section 6 we study the
wave front propagation properties of the solution $u(t,(x,k))$.

\section{Set up}

We shall first recall some basic facts in \cite{[F-H molecular
motors]} (also see \cite{[FW fish]}). Let us work with a fixed shape
of $D$.

First of all we need to construct a graph $\Gm$ related to the
domain $D$ (see Fig.1). For $x_0\in \R$, let $C(x_0)=\{(x,z)\in D:
x=x_0\}$ be the intersection of the domain $D$ with the line
$\{x=x_0\}$. The set $C(x_0)$ may have several connected components.
We identify all points in each connected component and the set thus
obtained, equipped with the natural topology, is homeomorphic to a
graph $\Gm$. We label the edges of this graph $\Gm$ by
$I_1,...,I_k,...$ (there might be infinitely many such edges).

We see that the structure of the graph $\Gm$ consists of many edges
(such as $I_1, I_3, I_5, I_7, I_9$,... in Fig.1) that form a long
line corresponding to the domain $D_0$ and many other short edges
(such as $I_2, I_4, I_6, I_8$,... in Fig.1) attached to the long
line in a random way. We will henceforth denote by $I_0$ the long
line corresponding to the domain $D_0$.

A point $y\in \Gm$ can be characterized by two coordinates: the
horizontal coordinate $x$, and the discrete coordinate $k$ being the
number of the edge $I_k$ in the graph $\Gm$ to which the point $y$
belongs. Let the identification mapping be $\iY: D \ra \Gm$. We note
that the second coordinate is not chosen in a unique way: for $y$
being an interior vertex $O_i$ of the graph $\Gm$ we can take $k$ to
be the number of any of the several edges meeting at the vertex
$O_i$.

The distance $\rho(y_1,y_2)$ between two points $y_1=(x_1, k)$ and
$y_2=(x_2, k)$ belonging to the same edge of the graph $\Gm$ is
defined as $\rho(y_1,y_2)=|x_1-x_2|$; for $y_1, y_2\in \Gm$
belonging to different edges of the graph it is defined as the
geodesic distance
$\rho(y_1,y_2)=\min(\rho(y_1,O_{j_1})+\rho(O_{j_1},
O_{j_2})+...+\rho(O_{j_l}, y_2))$, where the minimum is taken over
all chains $y_1\leftrightarrow O_{j_1}\leftrightarrow
O_{j_2}\leftrightarrow...\leftrightarrow O_{j_l}\leftrightarrow y_2$
of vertices $O_{j_i}$ connecting the points $y_1$ and $y_2$.

For an edge $I_k=\{(x,k): A_k\leq x\leq B_k\}$ we consider the
"tube" $U_k=\iiY^{-1}(I_k)\cap \{A_k\leq x\leq B_k\}$ in $D$. The
"tube" $U_k$ can be characterized by the interval $x\in [A_k,B_k]$
and the "height functions" $h_k^{\pm}(x)$: $U_k=\{(x,z): A_k\leq
x\leq B_k , h_k^-(x)\leq z \leq h_k^+(x)\}$. For $x\in [A_k, B_k]$,
we denote the set $C_k(x)$ to be the connected component of $C(x)$
that corresponds to the "tube" $U_k$: $C_k(x)=\{x\}\times [h_k^-,
h_k^{+}]$. Let $l_k(x)=h_k^+(x)-h_k^-(x)\geq 0$ for all $x\in \R$.
We notice that each $h_k^{\pm}(x) , l_k(x)$, etc. is smooth.

The vertices $O_j$ correspond to the connected components containing
points $(x,z)\in \pt D$ with $n_2(x,z)=0$. There are two types of
vertices: the interior vertices (in Fig.1 they are $O_1, O_3, O_5,
O_7$) are the intersection of three edges; the exterior vertices (in
Fig.1 they are $O_2, O_4, O_6, O_8$) are the endpoints of only one
edge.

Using the ideas in \cite{[FW fish]} with a little modification we
can establish the weak convergence of the process
$Y_t^\ve=\iiY(\boldX_t^\ve)$ (which is not Markov in general) as
$\ve \da 0$ in the space $\contfunc_{[0,T]}(\Gm)$ to a certain
Markov process $Y_t$ on $\Gm$. A sketch of the proof of this fact is
in \cite[Section 2]{[F-H molecular motors]}.

The process $Y_t$ is a diffusion process on $\Gm$ with a generator
$A$ and the domain of definition $D(A)$. We are going now to define
the operator $A$ and its domain of definition $D(A)$.

For each edge $I_k$ we define an operator $\overline{L}_k$:

$$\overline{L}_k u(x)=\dfrac{1}{2l_k(x)}\dfrac{d}{dx}
\left(l_k(x)\dfrac{du}{dx}\right)+\overline{V}_k(x)\dfrac{du}{dx}  \
, \ A_k \leq x \leq B_k \ . $$

Here
$$\overline{V}_k(x)=\dfrac{1}{l_k(x)}\int_{h_k^-(x)}^{h_k^+(x)}V(x,z)dz$$
is the average of the velocity field $V(x,z)$ on the connected
component $C_k(x)$, with respect to Lebesgue measure in
$z$-direction. At places where $l_k=0$, the above expression for
$\overline{V}_k(x)$ is understood as a limit as $l_k\ra 0$:
$$\overline{V}_k(x)=\lim\li_{y\ra x}\dfrac{1}{l_k(y)}\int_{h_k^-(y)}^{h_k^+(y)}V(y,z)dz \ .$$

We will assume throughout this paper the following.

\textbf{Assumption 7.} The function $\Vbar_k(x)=0$.

We notice that this is a bit different from the corresponding one in
\cite[Assumption 6]{[F-H molecular motors]}. We point out that the
vanishing mean drift assumption is crucial for the method of our
analysis to work.

Thus under our Assumption 7 we have $$\overline{L}_k
u(x)=\play{\dfrac{1}{2l_k(x)}\dfrac{d}{dx}\left(l_k(x)\dfrac{du}{dx}\right)}
\ .$$

The operator $\overline{L}_k$ can be represented as a generalized
second order differential operator (see \cite{[Feller]})

$$\overline{L}_k u(x)=D_{m_k}D_{p_k} f(x)  \ ,$$
where, for an increasing function $h$, the derivative $D_h$ is
defined by $D_h g(x)=\lim\li_{\dt \da
0}\dfrac{g(x+\dt)-g(x)}{h(x+\dt)-h(x)}$, and

$$p_k(x)=\int \dfrac{dx}{l_k(x)}$$ is the scale function, $$m_k(x)=2\int
l_k(x)dx$$ is the speed measure.

The operator $A$ is acting on functions $f$ on the graph $\Gm$: for
$y=(x,k)$ being an interior point of the edge $I_k$ we take
$Af(y)=\overline{L}_k f(x,k)$.

The domain of definition $D(A)$ of the operator $A$ consists of such
functions $f$ satisfying the following properties.

$\bullet$ The function $f$ is a continuous function that is twice
continuously differentiable in $x$ in the interior part of every
edge $I_k$;

$\bullet$ There exist finite limits $\lim\li_{y\ra O_i}A f(y)$
(which are taken as the value of the function $Af$ at the point
$O_i$);

$\bullet$ There exist finite one-sided limits $\lim\li_{x\ra
x_i}D_{p_k}f(x,k)$ along every edge ending at $O_i=(x_i,k)$ and
 they satisfy the gluing conditions

$$\sum\li_{j=1}^{N_i}(\pm)\lim\li_{x\ra x_i}D_{p_{k_j}}f(x,k_j)=0 \ , \eqno(2.1)$$
where the sign "$+$" is taken if the values of $x$ for points
$(x,k_j)\in I_{k_j}$ are $\geq x_i$ and "$-$" otherwise. Here
$N_i=1$ (when $O_i$ is an exterior vertex) or $3$ (when $O_i$ is an
interior vertex).

For an exterior vertex $O_i=(x_i,k)$ with only one edge $I_k$
attached to it the condition (2.1) is just $\lim\li_{x\ra
x_i}D_{p_k}f(x,k)=0$. Such a boundary condition can also be
expressed in terms of the usual derivatives $\dfrac{d}{dx}$ instead
of $D_{p_k}$. It is $\lim\li_{x\ra x_i}l_k(x)\dfrac{\pt f}{\pt
x}(x,k)=0$. We remark that we are in dimension 2 so that these
exterior vertices are accessible, and the boundary condition can be
understood as a kind of (not very standard) instantaneous
reflection. In dimension 3 or higher these endpoints do not need a
boundary condition, they are just inaccessible. For an interior
vertex the gluing condition (2.1) can be written with the
derivatives $\dfrac{d}{dx}$ instead of $D_{p_k}$. For $k$ being one
of the $k_j$ we define $\al_{ik}=\lim\li_{x\ra x_i}l_k(x)$ (for each
edge $I_k$ the limit is a one-sided one). Then the condition (2.1)
can be written as

$$\sum\li_{j=1}^3 (\pm)\al_{i,k_j}\cdot\lim\li_{x\ra x_i}\dfrac{df(x,k_j)}{dx}=0 \ . \eqno(2.2)$$

It can be shown as in \cite[Section 2]{[FW fish]} that the process
$Y_t$ exists as a continuous strong Markov process on $\Gm$.

We fix the shape of $D$. For every $\ve>0$, every $\boldx=(x,z)\in
D$ and every $T\in (0,\infty)$ let us consider the distribution
$\mu_{\boldx}^\ve$ of the trajectory $Y_t^\ve=\iiY(\boldX_t^\ve)$
starting from a point $\boldX_0^\ve=\boldx$ in the space
$\contfunc_{[0,T]}(\Gm)$ of continuous functions on the interval
$[0,T]$ with values in $\Gm$: the probability measure defined for
every Borel subset $B\subseteq \contfunc_{[0,T]}(\Gm)$ as
$\mu_{\boldx}^\ve(B)=\Prob_{\boldX_0^\ve=\boldx}^W(Y_{\bullet}^\ve\in
B)$. Similarly, for every $y\in \Gm$ and $T>0$ let $\mu^0_y$ be the
distribution of the process $Y_t$ in the same space:
$\mu^0_y(B)=\Prob^W_y(Y_{\bullet}\in B)$. The following theorem is
basic for our analysis.

\

\textbf{Theorem 2.1.} \textit{For every $\boldx\in D$ and every
$T>0$ the distribution $\mu_{\boldx}^\ve$ converges weakly to
$\mu_{\iiY(\boldx)}^0$ as $\ve\da 0$.}

\

In other words we have $$\E_{\boldX_0^\ve=\boldx}^W
F(Y^\ve_\bullet)\ra \E_{\iiY(\boldx)}^W F(Y_{\bullet})$$ for every
bounded continuous functional $F$ on the space
$\contfunc_{[0,T]}(\Gm)$.

The proof of this theorem follows from \cite{[FW fish]} and there is
a sketch in \cite[Section 3]{[F-H molecular motors]}. We omit
duplicating the details here.

\section{Convergence of $u^\ve$ to $u$}

We recall that our definition of the generalized solutions to (1.2)
and (1.4) are the solutions of the integral equations (1.3) and
(1.5), respectively. That is to say, we have

$$u^\ve(t,(x,z))=\E_{(x,z)}^W\left[\exp\left(\int_0^t
c(u^\ve(t-s, \boldX_s^\ve))ds\right)g(X_t^\ve)\right]$$ and
$$u(t,(x,k))=\E_{(x,k)}^W\left[\exp\left(\int_0^t c(u(t-s, Y_s))ds\right)
g(X_t)\right] \ .$$

\

\textbf{Theorem 3.1.} \textit{There exist unique bounded measurable
generalized solutions $u^\ve(t,(x,z))$, $t>0, (x,z)\in D$ and
$u(t,(x,k))$, $t>0, (x,k)\in \Gm$ for (1.2) and (1.4), respectively.
These solutions are continuous for all $t\geq 0$.}

\

\textbf{Proof.} We take (1.2) as an example. The proof for (1.4) is
exactly the same. We shall prove the existence and regularity by
using a contraction mapping principle (compare with \cite[\S 5.3]{[F
red book]}). To this end we consider the Banach space $B_T$ of
bounded measurable functions on $[0,T]\times [D]$ with norm
$\|v\|=\sup\li_{(t,x)\in [0,T]\times [D]}|v(t,x)|$. Consider in
$B_T$ the following operator

$$F[v]=F[v](t,(x,z))=\E_{(x,z)}^W\left[\exp\left(\int_0^t c(v(t-s, \boldX_s^\ve))ds\right)g(X_t^\ve)\right] \ , \ v\in B_T \ .$$

It is then checked that we have, for $0\leq t \leq T_0\leq T$, that

$$\begin{array}{l}
|F[u]-F[v]|(t, (x,z))\\
\play{=\left|\E_{(x,z)}^W\left[\exp\left(\int_0^t
c(u(t-s,\boldX_s^\ve))ds\right)g(X_t^\ve)\right]-\E_{(x,z)}^W\left[\exp\left(\int_0^t
c(v(t-s,\boldX_s^\ve)ds\right)g(X_t^\ve)\right]\right|}
\\
\play{\leq \|g\|\left|\E_{(x,z)}^W\left[\exp\left(\int_0^t
c(u(t-s,\boldX_s^\ve))ds\right)-\exp\left(\int_0^t
c(v(t-s,\boldX_s^\ve))ds\right)\right]\right|}
\\
\play{\leq \|g\|\exp(c(0)t)\cdot\text{Lip}(c)\cdot t \cdot \|u-v\|}
\\=C(T_0)\|u-v\|
\end{array}$$ for $u,v\in B_t$ and $C(T_0)=\|g\|\exp(c(0)T_0)\text{Lip}(c)T_0$. This
guarantees that $F$ is a contraction provided that
$T_0<\dfrac{1}{\|g\|\text{Lip}(c)\exp(c(0)T)}$. By contraction
mapping theorem we have existence and uniqueness of generalized
solution in the space of bounded measurable functions to the problem
(1.2) on the interval $[0, T_0]$. Since the solution $|u^\ve(T_0,
(x,z))|\leq \|g\|\exp(c(0)T_0)$ we can use the same $T_0$ and work
with intervals $[T_0, 2T_0]$, ..., up to $[(k-1)T_0, kT_0]$ ($k \geq
1$) provided that $kT_0<T$. This gives existence "in the large" for
a unique generalized solution $u^\ve(t, (x,z))$ for (1.2) in the
space of bounded measurable functions.

The continuity of the solution $u^\ve(t, (x,z))$ in the variables
$t$ and $(x,z)$ is provided by (1) The continuity of $g$; (2) The
Lipschitz continuity of $c(u)$; (3) The continuity and continuous
dependence of the process $\boldX_t^\ve$ on $\boldX_0^\ve=(x,z)$.
$\square$

\

We shall then show the approximation of the generalized solution
$u^\ve(t,(x,z))$ as $\ve$ is small by the generalized solution $u(t,
\iiY((x,z)))$. We prove this via a sequence of auxiliary results.

\

\textbf{Lemma 3.1.} \textit{We have}

$$\lim\li_{\ve \da 0}\max\li_{(x,z_1), (x,z_2)\in \iiY^{-1}((x,k))}
\max\li_{0\leq t \leq T}|u^\ve(t,(x,z_1))-u^\ve(t,(x,z_2))|=0 \ .
\eqno(3.1)$$

\

\textbf{Proof.} We consider the stopping time
$\tau^\ve=\tau^\ve((x,z_2), z_1)=\inf\{t>0: X_0^\ve=x, Z_0^\ve=z_2,
Z_t^\ve=z_1\}$. By strong Markov property of the process
$\boldX_t^\ve$ we have

$$\begin{array}{l}
u^\ve(t,(x,z_2))
\\
=\play{\E_{(x,z_2)}^W\left[\exp\left(\int_0^t c(u^\ve(t-s, X_s^\ve,
Z_s^\ve))ds\right)g(X_t^\ve)\right]}
\\
=\play{\E_{(x,z_2)}^W\left[\exp\left(\int_0^{\tau^\ve} c(u^\ve(t-s,
X_s^\ve, Z_s^\ve))ds\right)u^\ve(t-\tau^\ve, (X_{\tau^\ve}^\ve,
z_1))\right]} \ .
\end{array}$$

Since the motion $Z_t^\ve$ is moving very fast as $\ve\da 0$ we see
that $\tau^\ve\ra 0$ almost surely as $\ve \da 0$. This immediately
implies the convergence. $\square$

\

Let $(x,k)\in \Gm$. We introduce a new function $$\ubar^\ve(t,
(x,k))=\dfrac{1}{|\iiY^{-1}((x,k))|}\play{\int_{\iiY^{-1}((x,k))}u^\ve(t,(x,z))dz}
\ ,$$ and we see from Lemma 3.1 that we have the following.

\

\textbf{Corollary 3.1.} \textit{We have}
$$\lim\li_{\ve \da 0}\max\li_{(x,z)\in
\iiY^{-1}((x,k))}\max\li_{0\leq t \leq T}|\ubar^\ve(t,
\iiY(x,z))-u^\ve(t, (x,z))|=0 \ . \eqno(3.2)$$

\

We are going now to prove that the function $\ubar^\ve(t,(x,k))$ has
a uniform in $\ve$ bounded first derivative in the variable $x$.

\

\textbf{Lemma 3.2.} \textit{We have an a-priori estimate}
$$\max\li_{0\leq t \leq T, (x,k)\in \Gm}
\left|\dfrac{\pt \ubar^\ve}{\pt x}(t, (x,k))\right|\leq C(T)=C
\eqno(3.3)$$ \textit{where} $C>0$ \textit{is independent of} $\ve$.

\

\textbf{Proof.} By (1.3) we have

$$u^\ve(t, (x,z))=\E_{(x,z)}^W\left[\exp\left(\int_0^t
c(u^\ve(t-s, \boldX_t^\ve))ds\right)g(X_t^\ve)\right] \ .
$$

Differentiating with respect to $x$ we have

$$\begin{array}{l}
\play{\dfrac{\pt u^\ve}{\pt x}=\E^W\left[\dfrac{\pt g}{\pt
x}(X_t^\ve)\dfrac{\pt X_t^\ve}{\pt x} \exp\left(\int_0^t
c(u^\ve(t-s, X_s^\ve, Z_s^\ve))ds\right)\right.}
\\
\play{\left. \ \ \ +g(X_t^\ve)\exp\left(\int_0^t c(u^\ve(t-s,
X_s^\ve, Z_s^\ve))\right)\dfrac{\pt}{\pt x}\left(\int_0^t
c(u^\ve(t-s, X_s^\ve, Z_s^\ve))ds\right)\right]} \ .
\end{array}$$

Note that $$\dfrac{\pt}{\pt x}\left(\int_0^t c(u^\ve(t-s, X_s^\ve,
Z_s^\ve))ds\right)=\int_0^t c'(u^\ve(t-s, X_s^\ve,
Z_s^\ve))\dfrac{\pt u^\ve}{\pt x} (t-s, X_s^\ve, Z_s^\ve)\dfrac{\pt
X_s^\ve}{\pt x}ds \ .$$

Therefore if we let $m^\ve(t)=\max\li_{(x,k)\in D, 0\leq s \leq
t}\left|\dfrac{\pt u^\ve}{\pt x}(s, (x,k))\right|^2$, we get

$$m^\ve(t)\leq \al(t)+\bt(t)\int_0^t m^\ve(s)ds \ ,$$
where $\al(t), \bt(t)$ are bounded with their bound depending on the
regularity of $g(x)$, $c(u)$ and the shape parameter $l(x)$, yet
independent of $\ve$. We then apply a Gronwall inequality to
conclude that we have an a-priori estimate
$$\max\li_{0\leq t \leq T, (x,z)\in D}
\left|\dfrac{\pt u^\ve}{\pt x}(t, (x,z))\right|\leq C(T)=C
\eqno(3.4)
$$ where $C>0$ is independent of $\ve$.

From the above estimate and taking into account the smoothness of
the shape parameter $l(x)$, we see that the a-priori estimate in the
statement of the Lemma holds. In fact, by the definition of
$\ubar^\ve(t, (x,k))$, we have

$$|\iiY^{-1}((x,k))|\cdot \ubar^\ve(t,
(x,k))=\play{\int_{\iiY^{-1}((x,k))}u^\ve(t,(x,z))dz} \ .
$$

Thus $$\dfrac{\pt}{\pt
x}(|\iiY^{-1}((x,k))|)\ubar^\ve(t,(x,k))+|\iiY^{-1}((x,k))|\dfrac{\pt
\ubar^\ve}{\pt x}(t,(x,k))=\dfrac{\pt}{\pt
x}\left(\int_{\iiY^{-1}((x,k))}u^\ve(t,(x,z))dz\right) \ .$$

We notice that

$$\dfrac{\pt}{\pt
x}\left(\int_{\iiY^{-1}((x,k))}u^\ve(t,(x,z))dz\right)=\int_{\iiY^{-1}((x,k))}\dfrac{\pt
u^\ve}{\pt x}(t, (x,z))dz+u^\ve(t, (x,
b(x)))b'(x)-u^\ve(t,(x,a(x)))a'(x) \ .$$

Here $\iiY^{-1}(x,k)=(a(x), b(x))$. This implies (3.3). $\square$

\

Making use of Theorem 2.1, Corollary 3.1 and Lemma 3.2 we can prove
the following.

\

\textbf{Theorem 3.2.} \textit{We have}

$$\lim\li_{\ve \da 0}\max\li_{0\leq t \leq T}
\max\li_{(x,z)\in D}|u^\ve(t,(x,z))-u(t, \iiY((x,z)))|=0 \ .
\eqno(3.5)$$

\

\textbf{Proof.} An outline of this proof is mentioned at the end of
\cite[Chapter 7]{[F green book]}. We fulfill the details here.

Let $(x,k)=\iiY((x,z))$. From (1.3) and (1.5) we see that

$$\begin{array}{l}
|u^\ve(t, (x,z))-u(t, (x,k))|
\\
\play{=\left|\E_{(x,z)}^W\left[\exp\left(\int_0^t c(u^\ve(t-s,
\boldX_s^\ve))ds\right)g(X_t^\ve)\right]-\E_{(x,k)}^W\left[\exp\left(\int_0^t
c(u(t-s, Y_s))ds\right)g(X_t)\right]\right|}
\\
\leq \|g\|\cdot[(I)+(II)+(III)] \ .
\end{array}$$

Here

$$(I)=\left|\E_{(x,z)}^W\exp\left(\int_0^t c(u^\ve(t-s, \boldX_s^\ve))ds\right)-
\E_{(x,z)}^W\exp\left(\int_0^t c(\ubar^\ve(t-s,
\iiY(\boldX_s^\ve))ds\right)\right| \ ,$$

$$(II)=\left|\E_{(x,z)}^W\exp\left(\int_0^t c(\ubar^\ve(t-s, \iiY(\boldX_s^\ve)))ds\right)-
\E_{(x,k)}^W\exp\left(\int_0^t c(\ubar^\ve(t-s, Y_s)ds\right)\right|
\ ,$$

$$(III)=\left|\E_{(x,k)}^W\exp\left(\int_0^t c(\ubar^\ve(t-s, Y_s))ds\right)-
\E_{(x,k)}^W\exp\left(\int_0^t c(u(t-s, Y_s)ds\right)\right| \ .$$

Thus we see that $(I)\ra 0$ as $\ve \da 0$ due to (3.2); $(II)\ra 0$
as $\ve\da 0$ due to the weak convergence of the processes
$Y_t^\ve=\iiY(\boldX_t^\ve)$ to $Y_t$ on $\Gm$ as $\ve \da 0$ in
$\contfunc_{[0,T]}(\Gm)$ (Theorem 2.1) and (3.3); $(III)$ can be
bounded by a constant multiple of $\play{\int_0^t \max\li_{(x,z)\in
D}|u^\ve(s, (x,z))-u(s,(x,k))|ds}$ plus a term going to $0$ as $\ve
\da 0$ (due to (3.2)). We then apply a standard technique via
Gronwall's inequality and we can conclude. $\square$

\section{Auxiliary results}

This section will be devoted to obtaining some auxiliary results
which will be used in Sections 5--6.

Let the random variable

$$T_r^s=\inf\{t\geq 0, Y_0=(s,0), Y_t\in I_0, X_t\leq r\} \eqno(4.1.1)$$
for $s\geq r\in \R$. Intuitively, $T_r^s$ is the first time that the
process $Y_t$, starting from $Y_0=(s,0)$, comes back to $I_0$ with
the value of its $x$-component $\leq r$. We recall that $I_0$ is the
long line in $\Gm$ corresponding to the domain $D_0$.

In the same way we define

$$T_r^s=\inf\{t\geq 0, Y_0=(s,0), Y_t\in I_0, X_t\geq r\} \eqno(4.1.2)$$
for $s\leq r\in \R$.

Let $\lb>0$. Let the function
$$u(x)=\E^W [e^{-\lb T^x_0}] \eqno(4.2)$$ for $x\in \R$. We remind
the reader of a small notational convention here. In this section
for convenience of notation we have a minus sign in front of the
stopping time $T^x_0$ in (4.2). In the Sections 5--6 we will drop
this minus sign and instead we will be mainly working with $\lb<0$.

Let us first consider the case when $x>0$. The function $u(x)$ is
the solution of the following Sturm-Liouville problem

$$Au-\lb u=0 \text{ on } \Gm \ ,  \ u\in D(A) \ ,  \ u(0)=1 \ , \ u(+\infty)=0 \ . \eqno(4.3)$$

To solve the above problem we shall first recall the basic theory of
Feller (\cite{[Feller]}). We follow here \cite{[Mandl]} and we also
refer the reader to \cite[Lemma 2.10]{[FHW cone]}. Without loss of
generality let us first work with some interval $I=[0,r]$ for some
$r>0$. We consider the eigenvalue problem associated with the
generalized second-order differential operator $D_mD_p$:
$$D_mD_p u(x)-\lb u(x)=0 \eqno(4.4)$$ on an interval $x\in
I=[0,r]$. Here $m=m(x)$ is the speed measure and $p=p(x)$ is the
scale function. The function $p(x)$ is a strictly increasing
continuous function on $(0,r)$ and the function $m(x)$ is a strictly
increasing function on $(0,r)$ continuous to the right. The
generalized derivatives are defined as
$$D_pf(x)=\lim\li_{y\ra x}\dfrac{f(y)-f(x)}{p(y)-p(x)} \ , \ D_mf(x)=\lim\li_{y\ra x}\dfrac{f(y)-f(x)}{m(y)-m(x)}$$
where $x\in (0,r)$ and $f$ is a real function defined in a
neighborhood of $x$. There are two basic solutions $u_+(x)$,
$u_-(x)$ of the equation (4.4) with $u_+(0)=u_-(r)=0$ and
$u_+(r)=u_-(0)=1$; the function $u_+(x)$ is increasing in $x$ and
$u_-(x)$ is decreasing in $x$; the derivatives $D_p u_+(x)$, $D_p
u_-(x)$ are increasing functions.

Moreover, an explicit representation of the functions $u_{\pm}(x)$
is available (\cite{[Mandl]}). We set, for $n=0,1,2,...$, $x\in
[0,r]$,

$$u^{(0)}(x)\equiv 1 \ , \ u^{(n+1)}(x)=\int_0^x\int_0^y u^{(n)}(s)dm(s)dp(y) \ .$$

Let $$u(x,\lb)=\sum\li_{n=0}^\infty \lb^n u^{(n)}(x) \ .
\eqno(4.5.0)$$ It could be justified that the above series
converges.

We can easily check that $$D_p u(x,\lb)=\sum\li_{n=0}^\infty
\lb^{n+1}\int_0^x u^{(n)}(y)dm(y) \ . \eqno(4.5.0')$$

Let

$$u_+(x,\lb)=u(x,\lb)\int_0^x u(y,\lb)^{-2}dp(y) \ ,$$

$$u_-(x,\lb)=u(x,\lb)\int_x^r u(y,\lb)^{-2}dp(y) \ .$$

Moreover, we can calculate the derivative

$$\dfrac{du_+}{dx}(x,\lb)=\dfrac{du}{dx}(x,\lb)\int_0^x u(y,\lb)^{-2}dp(y)
+\dfrac{1}{u(x,\lb)}\dfrac{dp(x)}{dx} \ ,$$

$$\dfrac{du_-}{dx}(x,\lb)=\dfrac{du}{dx}(x,\lb)\int_x^r u(y,\lb)^{-2}dp(y)
-\dfrac{1}{u(x,\lb)}\dfrac{dp(x)}{dx} \ .$$

Then we have $$u_+(x)=\dfrac{u_+(x,\lb)}{u_+(r,\lb)} \ , \
u_-(x)=\dfrac{u_-(x,\lb)}{u_-(0,\lb)} , \eqno(4.5.1)$$ and

$$\begin{array}{l}
\play{\dfrac{du_+}{dx}(x)=\dfrac{\play{\dfrac{du}{dx}(x,\lb)\int_0^x
u(y,\lb)^{-2}dp(y)
+\dfrac{1}{u(x,\lb)}\dfrac{dp(x)}{dx}}}{\play{u(r,\lb)\int_0^r
u(y,\lb)^{-2}dp(y)}}  \ , }
\\
\play{\dfrac{du_-}{dx}(x)=\dfrac{\play{\dfrac{du}{dx}(x,\lb)\int_x^r
u(y,\lb)^{-2}dp(y)
-\dfrac{1}{u(x,\lb)}\dfrac{dp(x)}{dx}}}{\play{u(0,\lb)\int_0^r
u(y,\lb)^{-2}dp(y)}} \ . }
\end{array} \eqno(4.5.2)$$ In terms of generalized derivatives we
see that the above is equivalent to

$$\begin{array}{l}
\play{D_pu_+(x)=\dfrac{\play{D_pu(x,\lb)\int_0^x u(y,\lb)^{-2}dp(y)
+\dfrac{1}{u(x,\lb)}}}{\play{u(r,\lb)\int_0^r u(y,\lb)^{-2}dp(y)}} \
, }
\\
\play{D_pu_-(x)=\dfrac{\play{D_pu(x,\lb)\int_x^r u(y,\lb)^{-2}dp(y)
-\dfrac{1}{u(x,\lb)}}}{\play{u(0,\lb)\int_0^r u(y,\lb)^{-2}dp(y)}} \
. }
\end{array} \eqno(4.5.2')$$

A general solution of (4.4) can be represented as a linear
combination
$$u(x)=c^+u_+(x)+c^-u_-(x) \ .$$ The constants $c^+$ and $c^-$ are
determined by boundary conditions to be specified.

In the case when $r<0$ situation is similar and we have to make
small changes accordingly. To be more precise, we can treat the
point $r$ as the point $0$ and the point $0$ as the point $r$. The
formulas (4.5.0) ((4.5.0$'$)), (4.5.1) and (4.5.2) ((4.5.2$'$)) have
to be changed accordingly. In the rest of this section we will be
mainly performing detailed steps in the calculations assuming $r>0$
and we will present corresponding results when $r<0$ without a
detailed calculation.

Let us come back to our problem (4.3). First of all we note that the
structure of the graph $\Gm$ consists of two types of edges: the
first type of edges are lined up together forming the edge $I_0$ and
we label them as $I_{2k-1}$, $k\in \N$; the second type of edges
correspond to the wings and we label them as $I_{2k}$, $k\in \N$.
These edges are labeled in a consecutive way (see Fig.1). Let the
projection of the second type of edges $I_{2k}$ onto the
$x$-direction be isomorphic to $[0,r_{2k}]$ for $r_{2k}>0$
 and $[r_{2k},0]$ for $r_{2k}<0$. Let the interval $I_{2k-1}$ be isomorphic to $[0, r_{2k-1}]$.
 We solve the problem $\overline{L}_ku_k-\lb u_k=0$ on each edge
$I_{2k-1}\cong [0,r_{2k-1}]$ (the first type), $I_{2k}\cong
[0,r_{2k}]$ for $r_{2k}>0$
 and $I_{2k+1}\cong [r_{2k},0]$ for $r_{2k}<0$ (the second type). We notice that in this case
when we represent the operator $\overline{L}_k$ as a generalized
second order derivative operator $\overline{L}_k=D_{m_k}D_{p_k}$ we
will have
$$dm_k(x)=2l_k(x)dx \eqno(4.5.3)$$ and $$dp_k(x)=\dfrac{1}{l_k(x)}dx \ . \eqno(4.5.4)$$

The general solution is represented as $u_k(x)=c^+_k u_+^k(x)+c^-_k
u_-^k(x)$. Here $u_+^k(x)$ and $u_-^k(x)$ are the two basic
solutions corresponding to the interval $I_k$ and we identify $x$
with some $x\in I_k$ (or its projection onto the $x$-axis, anyway).
We note that they are random solutions. The constants $c^+_k$ and
$c^-_k$ are to be determined. We shall seek for a solution
$u(x)=u_{2k-1}(x)$ whenever $(x,0)\in \Gm$. Thus we have

$$\left\{\begin{array}{l}
D_{p_{2k}}u_{2k}(r_{2k})=0 \ ,
\\
u_{2k-1}(r_{2k-1})=u_{2k+1}(0)=u_{2k}(0) \ ,
\\
\al_{2k-1}\dfrac{du_{2k-1}}{dx}(r_{2k-1})
=\text{sign}(r_{2k})\gm_{2k}\dfrac{du_{2k}}{dx}(0)+\bt_{2k+1}\dfrac{du_{2k+1}}{dx}(0)
\ ,
\\
u_1(0)=1  \ ,
\\
\lim\li_{k\ra\infty}u_{2k+1}(r_{2k+1})=0 \ .
\end{array}\right. \eqno(4.6)$$

In the above $\al_{2k-1}, \bt_{2k+1}, \gm_{2k}$ are the
corresponding cross-section width of the channel at the junctions.
We have $\al_{2k-1}-\bt_{2k+1}=\text{sign}(r_{2k})\gm_{2k}$.

\

\textbf{Lemma 4.1.} \textit{We have}
$$\begin{pmatrix}c^+_{2k-1}\\c^-_{2k-1}\end{pmatrix}=M_k
\begin{pmatrix}c^+_{2k+1}\\c^-_{2k+1}\end{pmatrix}  \ . \eqno(4.7)$$
\textit{Here}

$$M_k=\begin{pmatrix}0&1\\ x_k&y_k\end{pmatrix}$$ \textit{with}

$$x_k=\dfrac{\bt_{2k+1}}{\al_{2k-1}}\dfrac{\dfrac{d u_+^{2k+1}}{dx}(0)}{\dfrac{d u_-^{2k-1}}{dx}(r_{2k-1})} \ , \eqno(4.8)$$

$$\begin{array}{l}
y_k=\play{-\dfrac{\dfrac{du_+^{2k-1}}{dx}(r_{2k-1})}{\dfrac{du_-^{2k-1}}{dx}(r_{2k-1})}}
\\
\play{\ \ \ +\dfrac{\gm_{2k}}{\al_{2k-1}}\left(
-\dfrac{D_{p_{2k}}u_-^{2k}(r_{2k})}{D_{p_{2k}}u_+^{2k}(r_{2k})}
\dfrac{\dfrac{du_+^{2k}}{dx}(0)}{\dfrac{du_-^{2k-1}}{dx}(r_{2k-1})}+
\dfrac{\dfrac{du_-^{2k}}{dx}(0)}{\dfrac{du_-^{2k-1}}{dx}(r_{2k-1})}\right)}
\\
\play{\ \ \
+\dfrac{\bt_{2k+1}}{\al_{2k-1}}\dfrac{\dfrac{du_-^{2k+1}}{dx}(0)}{\dfrac{du_-^{2k-1}}{dx}(r_{2k-1})}}
\ ,
\end{array} \eqno(4.9.1)$$ \textit{if $r_{2k}>0$; and}

$$\begin{array}{l}
y_k=\play{-\dfrac{\dfrac{du_+^{2k-1}}{dx}(r_{2k-1})}{\dfrac{du_-^{2k-1}}{dx}(r_{2k-1})}}
\\
\play{\ \ \ -\dfrac{\gm_{2k}}{\al_{2k-1}}\left(
\dfrac{\dfrac{du_+^{2k}}{dx}(0)}{\dfrac{du_-^{2k-1}}{dx}(r_{2k-1})}-\dfrac{D_{p_{2k}}u_+^{2k}(r_{2k})}{D_{p_{2k}}u_-^{2k}(r_{2k})}
\dfrac{\dfrac{du_-^{2k}}{dx}(0)}{\dfrac{du_-^{2k-1}}{dx}(r_{2k-1})}\right)}
\\
\play{\ \ \
+\dfrac{\bt_{2k+1}}{\al_{2k-1}}\dfrac{\dfrac{du_-^{2k+1}}{dx}(0)}{\dfrac{du_-^{2k-1}}{dx}(r_{2k-1})}}
\ ,
\end{array} \eqno(4.9.2)$$ \textit{if $r_{2k}<0$.}

\textbf{Proof.} The first three equalities in (4.6) will give us

$$\begin{pmatrix}c^+_{2k-1}\\c^-_{2k-1}\end{pmatrix}=M_k \begin{pmatrix}c^+_{2k+1}\\c^-_{2k+1}\end{pmatrix}$$
where we can calculate the random matrix $M_k$. Let us first
consider the case when $r_{2k}>0$. From the first equation of (4.6)
we see that we have

$$c_{2k}^+D_{p_{2k}}u^{2k}_+(r_{2k})+c_{2k}^-D_{p_{2k}}u^{2k}_-(r_{2k})=0 \ . \eqno(4.10)$$

We have
$u_{2k-1}(r_{2k-1})=c_{2k-1}^+u_+^{2k-1}(r_{2k-1})=c_{2k-1}^+$,
$u_{2k+1}(0)=c_{2k+1}^- u_-^{2k+1}(0)=c_{2k+1}^-$,
$u_{2k}(0)=c_{2k}^- u_-^{2k}(0)=c_{2k}^-$. So from the second
equality of (4.6) we get

$$c_{2k-1}^+=c_{2k+1}^-=c_{2k}^- \ . \eqno(4.11)$$

The third equality in (4.6) gives us

$$\begin{array}{l}
\play{\al_{2k-1}\left(c_{2k-1}^+\dfrac{du_+^{2k-1}}{dx}(r_{2k-1})+c_{2k-1}^-\dfrac{du_-^{2k-1}}{dx}(r_{2k-1})\right)}
\\
\play{=\gm_{2k}\left(c_{2k}^+\dfrac{du_+^{2k}}{dx}(0)+c_{2k}^-\dfrac{du_-^{2k}}{dx}(0)\right)
+\bt_{2k+1}\left(c_{2k+1}^+\dfrac{du_+^{2k+1}}{dx}(0)+c_{2k+1}^-\dfrac{du_-^{2k+1}}{dx}(0)\right)
\ .}
\end{array} \eqno(4.12)$$

From (4.10), (4.11) and (4.12) we can conclude that

$$M_k=\begin{pmatrix}0&1\\ x_k&y_k\end{pmatrix} \ .$$

Here the random variables $x_k$ and $y_k$ are defined by (4.8) and
(4.9.1).

In the case when $r_{2k}<0$ we just have to replace in (4.11) the
coefficient $c_{2k}^-$ by $c_{2k}^+$ and we have to change the sign
in front of $\gm_{2k}$ in (4.12) to minus. We thus get (4.9.2).
$\square$

\

Let $L_{k}=\sum\li_{j=1}^k r_{2j-1}$. Due to shift invariance
(stationarity) we see that without loss of generality we can assume
that $L_1$ and $L$ have the same distribution.

\

\textbf{Lemma 4.2.} \textit{The random variables}
$\dfrac{c_{2k-1}^+}{c_{2k-1}^-}$ , $k\in \N$ \textit{are identically
distributed. }

\

\textbf{Proof.} We have, by strong Markov property of the process
$Y_t$ on $\Gm$, that for $\lb>0$,

$$\E^W[e^{-\lb T^{L_{k+1}}_0}]=\E^W[e^{-\lb T^{L_{k+1}}_{L_{k}}}]\E^W[e^{-\lb T^{L_k}_0}] \ .$$

By stationarity we see that $\E^W[e^{-\lb T^{L_{k+1}}_{L_{k}}}]$ has
the same distribution as $\E^W[e^{-\lb T^{L_1}_0}]$. Since
$$\E^W[e^{-\lb T^{L_{k+1}}_0}]=u(L_{k+1})=c_{2k-1}^+$$ and
$$\E^W[e^{-\lb T^{L_k}_0}]=u(L_k)=c_{2k-1}^- \ , $$ we see that the
distribution of $\dfrac{c_{2k-1}^+}{c_{2k-1}^-}$ is independent of
$k\in \N$. $\square$

\

\textbf{Lemma 4.3.} \textit{For any $k\in \N$ we have $0\leq
\dfrac{c_{2k+1}^+}{c_{2k+1}^-}\leq 1$ almost surely with respect to
$\Prob$.}

\

\textbf{Proof.} This is because we have $\E^W[e^{-\lb T_0^{L_1}}]\in
[0,1]$ and Lemma 4.2. $\square$

\

It is convenient to introduce the notation
$$S_{k}=\text{sign}(r_k)\int_0^{r_k}u_k(y,\lb)^{-2}\dfrac{dy}{l_k(y)} \ .$$ Here
$u_k(x,\lb)$ is given by the formula (4.5.0) specified in the
interval $I_k$.

\

\textbf{Lemma 4.4.} \textit{For any $k\in \N$ we have
$$x_k=-\dfrac{u_{2k-1}(r_{2k-1},\lb)S_{2k-1}}{u_{2k+1}(r_{2k+1},\lb)S_{2k+1}}
\ ;
$$ $y_k>0$ almost surely with respect to $\Prob$.}

\

\textbf{Proof.} By Lemma 4.1 we have
$$\dfrac{c_{2k-1}^+}{c_{2k-1}^-}=\dfrac{1}{x_k\dfrac{c_{2k+1}^+}{c_{2k+1}^-}+y_k} \ .$$

Making use of (4.5.2) and (4.5.3) we see that
$\play{x_k=-\dfrac{u_{2k-1}(r_{2k-1},\lb)S_{2k-1}}{u_{2k+1}(r_{2k+1},\lb)S_{2k+1}}}$.
Thus $x_k<0$ and since we have Lemma 4.2 and Lemma 4.3 we see that
$y_k>0$, almost surely with respect to $\Prob$. $\square$

\

\textbf{Lemma 4.5.} \textit{ We have  $\dfrac{1}{y_k}\leq c_1^+ \leq
1$ and thus $y_k\geq 1$ almost surely with respect to $\Prob$.}

\

\textbf{Proof.} This is a simple consequence of Lemma 4.1, Lemma 4.3
and Lemma 4.4. $\square$

\

\textbf{Lemma 4.6.} \textit{We have} $$\begin{array}{l} y_k
\play{=\dfrac{1}{u_{2k-1}(r_{2k-1},\lb)}+\al_{2k-1}\dfrac{du_{2k-1}}{dx}(r_{2k-1},\lb)S_{2k-1}}
\\ \ \ \ +u_{2k-1}(r_{2k-1},\lb)\dfrac{S_{2k-1}}{S_{2k}}
\dfrac{u_{2k}(r_{2k},\lb)D_{p_{2k}}u_{2k}(r_{2k},\lb)S_{2k}}
{u_{2k}(r_{2k},\lb)D_{p_{2k}}u_{2k}(r_{2k},\lb)S_{2k}+1}+u_{2k-1}(r_{2k-1},\lb)\dfrac{S_{2k-1}}{S_{2k+1}}
\ , \end{array} \eqno(4.13.1)$$ \textit{when $r_{2k}>0$; and}
$$\begin{array}{l} y_k
\play{=\dfrac{1}{u_{2k-1}(r_{2k-1},\lb)}+\al_{2k-1}\dfrac{du_{2k-1}}{dx}(r_{2k-1},\lb)S_{2k-1}}
\\ \ \ \ +\gm_{2k}\dfrac{u_{2k-1}(r_{2k-1},\lb)}{u_{2k}(0,\lb)}\dfrac{du_{2k}}{dx}(0,\lb)S_{2k-1}+u_{2k-1}(r_{2k-1},\lb)\dfrac{S_{2k-1}}{S_{2k+1}}
\ , \end{array} \eqno(4.13.2)$$ \textit{when $r_{2k}<0$.}

\

\textbf{Proof.} Let us first consider the case when $r_{2k}>0$. By
(4.9.1) we can write
$$y_k=(i)+\dfrac{\gm_{2k}}{\al_{2k-1}}((ii)(iii)+(iv))+\dfrac{\bt_{2k+1}}{\al_{2k-1}}(v) \ .$$

Here
$$(i)=-\dfrac{\dfrac{du^{2k-1}_+}{dx}(r_{2k-1})}{\dfrac{du_-^{2k-1}}{dx}(r_{2k-1})} \
, (ii)=-\dfrac{D_{p_{2k}}u^{2k}_-
(r_{2k})}{D_{p_{2k}}u_+^{2k}(r_{2k})} \ ,
(iii)=\dfrac{\dfrac{du^{2k}_+}{dx}(0)}{\dfrac{du_-^{2k-1}}{dx}(r_{2k-1})}
\ ,
$$ $$(iv)=\dfrac{\dfrac{du^{2k}_-}{dx}(0)}{\dfrac{du_-^{2k-1}}{dx}(r_{2k-1})}
\ ,
(v)=\dfrac{\dfrac{du^{2k+1}_-}{dx}(0)}{\dfrac{du_-^{2k-1}}{dx}(r_{2k-1})}
\ .$$

By making use of (4.5.2) and $(4.5.2')$ as well as $(4.5.0')$ it is
straightforward to calculate that

$$(i)=\al_{2k-1}\dfrac{du_{2k-1}}{dx}(r_{2k-1},\lb)S_{2k-1}+\dfrac{1}{u_{2k-1}(r_{2k-1},\lb)} \ ,$$

$$(ii)=\dfrac{u_{2k}(r_{2k},\lb)}{u_{2k}(r_{2k},\lb)D_{p_{2k}}u_{2k}(r_{2k},\lb)S_{2k}+1} \ , $$

$$(iii)=-\dfrac{\al_{2k-1}}{\gm_{2k}}\dfrac{u_{2k-1}(r_{2k-1}, \lb)S_{2k-1}}{u_{2k}(r_{2k},\lb)S_{2k}} \ ,$$

$$(iv)=\dfrac{\al_{2k-1}}{\gm_{2k}}
u_{2k-1}(r_{2k-1},\lb)\dfrac{S_{2k-1}}{S_{2k}} \ ,$$

$$(v)=\dfrac{\al_{2k-1}}{\bt_{2k+1}}u_{2k-1}(r_{2k-1},\lb)\dfrac{S_{2k-1}}{S_{2k+1}}
\ .$$

So we get

$$\begin{array}{l}
y_{k}
\\
=\al_{2k-1}\dfrac{du_{2k-1}}{dx}(r_{2k-1},\lb)S_{2k-1}+\dfrac{1}{u_{2k-1}(r_{2k-1},\lb)}
\\
\ \
-\play{\dfrac{1}{u_{2k}(r_{2k},\lb)D_{p_{2k}}u_{2k}(r_{2k},\lb)S_{2k}+1}u_{2k-1}(r_{2k-1},
\lb)\dfrac{S_{2k-1}}{S_{2k}}}
\\
\
 \ +u_{2k-1}(r_{2k-1},\lb)\dfrac{S_{2k-1}}{S_{2k}}+u_{2k-1}(r_{2k-1},\lb)\dfrac{S_{2k-1}}{S_{2k+1}} \ .
\end{array}$$

Thus we get (4.13.1). The equality (4.13.2) is obtained in a similar
way. $\square$

\

Making use of Lemma 4.6 and basic calculations (4.5.0)--(4.5.4), as
well as our Assumptions 2 and 6, we see that we have the following.

\

\textbf{Corollary 4.1.} $$\E \ln y_k<\infty \ .$$

\

Combining Lemmas 4.5 and 4.6 we colculde that $\E|\ln c_1^+|\leq
\E\ln y_k<\infty$. Thus we have the following.

\

\textbf{Theorem 4.1.} $$\E[|\ln \E^W[e^{-\lb T_0^{L}}]|]<\infty \ ;
\ \E[|\ln \E^W[e^{-\lb T_0^{-L}}]|]<\infty \ .$$

\

In the above theorem the second inequality is estimated in a similar
fashion as the first one.

\

\textbf{Lemma 4.7.} \textit{We have $\Prob^W(T_0^L<\infty)=1$ and
$\Prob^W(T_0^{-L}<\infty)=1$ almost surely with respect to $\Prob$.}

\

\textbf{Proof.} We take $T_0^L$ as an example. The case for
$T^{-L}_0$ is similar. Let $v((x,k);
A)=\Prob^W(T_{[0,A]}^{(x,k)}<\infty)$ where $(x,k)\in \Gm$. Here
$T_{[0,A]}^{(x,k)}$ is the first time the process $Y_t$, starting
from $Y_0=(x,k)\in \Gm$, hits the point $(0,k=0)\in \Gm$ or $(A,
k=0)\in \Gm$. We set $v_k(x)=v((x,k); A)$ where we identify
$(x,k)\in \Gm$ with some $x\in I_k$ (or projection of $I_k$ onto the
$x$-axis, anyway). Then we have

$$\left\{\begin{array}{l}
\play{\dfrac{1}{2}\dfrac{1}{l_k(x)}\dfrac{d}{dx}\left(l_k(x)\dfrac{d
v_k}{dx}\right)=0 \ ,}
\\
\play{\lim\li_{x\ra r_{2k}}l_{2k}(x)\dfrac{dv_{2k}}{dx}(x)=0 \ ,}
\\
\play{\al_{2k-1}\dfrac{dv_{2k-1}}{dx}(r_{2k-1})=\text{sign}(r_{2k})\gm_{2k}\dfrac{dv_{2k}}{dx}(0)+\bt_{2k+1}\dfrac{dv_{2k+1}}{dx}(0)
\ ,}
\\
v_{2k-1}(r_{2k-1})=v_{2k+1}(0)=v_{2k}(0) \ ,
\\
\play{v_1(0)=1 \ ,}
\\
v((A,0); A)=0 \ .
\end{array}\right.$$

This gives

$$v_{2k}(x)=v_{2k}(0)=v_{2k-1}(r_{2k-1})=v_{2k+1}(0) \ ,$$

$$\dfrac{dv_{2k-1}}{dx}(x)=\dfrac{C_{2k-1}}{l_{2k-1}(x)} \ ,$$

$$C_{2k+1}=C_{2k-1} \ .$$

Thus

$$v_{2k+1}(x)=D_{2k+1}+C\int_0^x \dfrac{dy}{l_0(y)} \ .$$

We see that $D_1=1$. And we have recursively that
$$D_{2k+1}-D_{2k-1}=C\int_0^{r_{2k-1}}\dfrac{dy}{l_0(y)} \ .$$ Thus
we see that
$$v((x,0); A)=1-\dfrac{\play{\int_0^x\dfrac{dy}{l_0(y)}}}{\play{\int_0^A\dfrac{dy}{l_0(y)}}}
\ .$$ Thus $\lim\li_{A\ra\infty}v((x,0); A)=1$ by Assumption 1
($\infty>\bar{l}_0\geq l_0(x)\geq l_0>0$). In particular,
$\Prob^W(T_0^L<\infty)=1$ almost surely with respect to $\Prob$.
$\square$

\

\textbf{Lemma 4.8.} \textit{We have $\E^W[T_0^L]=+\infty$ and
$\E^W[T_0^{-L}]=+\infty$ almost surely with respect to $\Prob$.}

\

\textbf{Proof.} We take $T_0^L$ as an example. The proof for
$T_0^{-L}$ is similar. We show by comparison. To this end we
construct the \textit{part} of the process $\boldX_t^\ve=(X_t^\ve,
Z_t^\ve)$ within the domain $D_0$ (compare with \cite[Section
4.2]{[F-H molecular motors]}). Let $\play{\vphi_t=\int_0^t
\1(\boldX_s^\ve\in D_0)ds}$ be an additive functional, which is
called \textit{the proper time} of the domain $D_0$. We introduce
the time $\bt_t$ inverse to $\vphi_t$ and continuous on the right.
Let $Y_t^{\ve,D_0}=\iiY(\boldX_{\bt_t}^\ve)$. One can show that as
$\ve\da0$ the weak convergence of $Y_t^{\ve, D_0}$ to $Y_t^{D_0}$.
The process $Y_t^{D_0}$ is described as a one-dimensional diffusion
process on $\R$ with gluing conditions (see \cite[Section 4.2]{[F-H
molecular motors]}). We have $T_0^L\geq T_0^{L, I_0}$ where $T_0^{L,
I_0}$ is the proportion of time of process $Y_t$ spent inside $I_0$.
We see that $T_0^{L, I_0}=T_0^{L, D_0}$ where $T_0^{L,
D_0}=\inf\{t\geq 0: Y_0^{D_0}=L, Y_t^{D_0}=0\}$. It is not hard to
prove, via an approximation similar as in \cite[Section 4.2]{[F-H
molecular motors]}, that $\E^W[T_0^{L, D_0}]=+\infty$. More
precisely, let
$$v(x;A)=-2\int_0^x \dfrac{dy}{l_0(y)}\int_0^yl_0(z)dz+
2\dfrac{\play{\int_0^A \dfrac{dy}{l_0(y)}\int_0^y
l_0(z)dz}}{\play{\int_0^A \dfrac{dy}{l_0(y)}}}\int_0^x
\dfrac{dy}{l_0(y)} \ .$$ Then we have $\E^W[T_0^{L,
D_0}]=\lim\li_{A\ra\infty}v(L; A)$. Since we can estimate
$$\dfrac{\play{\int_0^A \dfrac{dy}{l_0(y)}\int_0^y
l_0(z)dz}}{\play{\int_0^A \dfrac{dy}{l_0(y)}}}\geq
\dfrac{A}{2}\dfrac{l_0^2}{\bar{l}_0}\ ,$$ and we have our Assumption
1, we see that $\E^W[T_0^L]=+\infty$. $\square$

\section{The Large deviation principle}

We are interested in describing the wave front propagation
corresponding to the solution $u(t,(x,k))$ of (1.5). To this end we
study the quenched large deviation principle for the random variable
$\dfrac{vt-X^{(vt,k)}(\kp t)}{\kp t}$. Here $v>0$, $\kp>0$ and
$X^{(vt,k)}(\kp t)$ is the first component of the process $Y_t=(X_t,
k_t)$ on $\Gm$ starting from a point $(vt, k)\in \Gm$. Here $k$ may
be $0$ or some other integer $\geq 1$ depending on the structure of
$\Gm$. This is in essence an adaptation of the arguments of
\cite{[Taleb]} and \cite{[Zeitouni LDP RWRE]}.

Here and below, for notational convenience we will use the symbol
$X^x(\kp t)$ to denote the process $X_t$ (which is the first
component of the process $Y_t=(X_t, k_t)$) starting from a point
$(x, k)$ on $\Gm$ with an arbitrary choice of $k$. The fact that the
large deviation results for the random variable
$\dfrac{vt-X^{(vt,k)}(\kp t)}{\kp t}$ are independent of the choice
of $k$ will be revealed in the proof of Theorem 5.2.

Let $\lb\in \R$ and we introduce $$q(r,s,\lb)=\E^W[e^{\lb
T_r^s}\1_{T_r^s<\infty}] \ . \eqno(5.1)$$

Recall that $L$ is the distance between two consecutive vertices
$O_i$ and $O_j$ at which there is an edge corresponding to a wing.
We see that $L$ is a random variable measurable with respect to the
filtration $\{\cF_s^t\}_{-\infty\leq t\leq s\leq \infty}$ generated
by the shape of $D$. For each fixed shape of $D$ the random
variables $T^L_0$ and $T^{-L}_0$ are well defined and they are
measurable with respect to the filtration generated by the Wiener
process $(W_t^1, W_t^2)$. Notice that by our Assumption 6 we have
$\infty>\overline{L}>L> \underline{L}>0$ where $\underline{L},
\overline{L}$ are constants.

\

\textbf{Lemma 5.1.} \textit{Suppose that $\lb\in \R$ is such that
$$\E\left(|\ln \E^W[e^{\lb T_0^L}\1_{T_0^L<\infty}]|\right)<\infty \
; \ \E\left(|\ln \E^W[e^{\lb
T_0^{-L}}\1_{T_0^{-L}<\infty}]|\right)<\infty \ .$$} \textit{Let
$\lb \in \R$ and $c<v$. Then almost surely with respect to $\Prob$
the limits}
$$\mu^+(\lb)\equiv\lim\li_{t\ra\infty}\dfrac{1}{(v-c)t}\ln \E^W[e^{\lb
T_{ct}^{vt}}\1_{T_{ct}^{vt}<\infty}]=\dfrac{1}{\E L}\E\left(\ln
\E^W[e^{\lb
T_0^L}\1_{T_0^L<\infty}]\right)$$
$$\mu^-(\lb)\equiv\lim\li_{t\ra\infty}\dfrac{1}{(v-c)t}\ln
\E^W[e^{\lb T_{vt}^{ct}}\1_{T_{vt}^{ct}<\infty}]=\dfrac{1}{\E
L}\E\left(\ln \E^W[e^{\lb T_0^{-L}}\1_{T_0^{-L}<\infty}]\right)$$
\textit{hold. The convergence is uniform with respect to $v$ and $c$
as $v$ and $c$ vary in a set that is bounded and $(v-c)$ is bounded
away from zero. Moreover, $\mu^\pm(\lb)$ is independent of $v$ and
$c$.}

\

\textbf{Proof.} Let us just work with $\mu^+(\lb)$. The proof of
this fact is essentially the same as that of \cite[Section 2,
Proposition 1]{[Nolen-Xin KPP random drift]} provided we make small
modifications. In fact, by the strong Markov property of the process
$Y_t$ on $\Gm$ it is easy to deduce that for $r<s<t$ we have $$\ln
q(r,t,\lb)=\ln q(r,s,\lb)+\ln q(s,t,\lb) \ .$$ Let there be located
$N(n)$ edges that correspond to the "wings" in the interval $x\in
[0, cn]$. We see that $\lim\li_{n\ra\infty}\dfrac{cn}{N(n)}=\E L$
holds $\Prob$--a.s.. On the other hand, we have, by the ergodic
theorem, that
$$\lim\li_{n\ra\infty}\dfrac{\ln q(0, cn ,\lb)}{N(n)}=\E\left(\ln
\E^W[e^{\lb T_0^L}\1_{T_0^L<\infty}]\right)$$ holds $\Prob$--a.s..
Thus we see that
$$\lim\li_{n\ra\infty}\dfrac{1}{cn}\ln q(0, cn ,\lb)
=\lim\li_{n\ra\infty}\dfrac{1}{\dfrac{cn}{N(n)}}\dfrac{\ln q(0, cn
,\lb)}{N(n)} =\dfrac{1}{\E L}\E\left(\ln \E^W[e^{\lb
T_0^L}\1_{T_0^L<\infty}]\right)$$ provided that
$$\E\left(|\ln \E^W[e^{\lb
T_0^L}\1_{T_0^L<\infty}]|\right)<\infty \ .$$ The rest of the
argument is the same as in \cite[Section 2, Proposition
1]{[Nolen-Xin KPP random drift]}. $\square$

\

We note that by Theorem 4.1 the requirements of Lemma 5.1 always
hold for $\lb< 0$.

Let $$\lb_c^{\pm}=\sup\{\lb\in \R, \mu^{\pm}(\lb)<\infty\} \ .$$ Our
Theorem 4.1 implies that $\lb_c^{\pm}\geq 0$.

\

\textbf{Lemma 5.2.} (Properties of the function $\mu^{\pm}(\lb)$)

(1)\textit{ $\mu^{\pm}(0)=0$;}

(2)\textit{ $\mu^{\pm}(\lb)<0$ for $\lb<0$;}

(3)\textit{ $\mu^{\pm}(\lb)\ra-\infty$ as $\lb \ra -\infty$;}

(4)\textit{ $\mu^{\pm}(\lb)=+\infty$ as $\lb>\lb_c^{\pm}$;}

(5)\textit{ $\mu^{\pm}(\lb)$ is convex for $\lb \in (-\infty,
\lb_c^{\pm})$;}

(6)\textit{ For $\lb<\lb_c^{\pm}$, $\mu^{\pm}(\lb)$ is differentiable and
$$(\mu^{\pm})'(\lb)=\E\left[\dfrac{\E^W[T_0^{\pm L} e^{\lb T_0^{\pm L}}\1_{T_0^{\pm L}<\infty}]}
{\E^W[e^{\lb T_0^{\pm L}}\1_{T_0^{\pm L}<\infty}]}\right] \ ; $$ In
particular, $$(\mu^\pm)'(0)=a_0^{\pm}=\E[\E^W[T_0^{\pm
L}\1_{T_0^{\pm L}<\infty}]]\in (0,\infty] \ ; $$}

(7)\textit{ $(\mu^\pm)'(\lb)$ is monotonically strictly increasing
for $\lb\in (-\infty, \lb_c)$; }

(8)\textit{ $a_0^+=a_0^-=+\infty$ and therefore
$\lb_c^+=\lb_c^-=0$.}

\

\textbf{Proof.} The proof of this lemma is the same as in
\cite[Lemma 2.2, Proposition 4.2]{[Nolen-Xin KPP random drift]}. The
last statement (8) follows from our Lemmas 4.7 and 4.8. $\square$

\

We define

$$I^\pm(a)\equiv\sup\li_{\lb \leq 0}(a\lb-\mu^\pm(\lb)) \ .$$

\

\textbf{Lemma 5.3.} (Properties of the function $I^{\pm}(a)$)

(1) \textit{$I^{\pm}(a)>0$ for $a\in (0,\infty)$;}

(2) \textit{$I^{\pm}(a)$ is convex and decreasing in $a$ for $a\in
(0,\infty)$;}

(3) \textit{$\lim\li_{a\ra 0+}I^{\pm}(a)=+\infty$ and $\lim\li_{a\ra
+\infty}I^\pm(a)=0$.}

\

\textbf{Proof.} The proof of this lemma is the same as in
\cite{[Nolen-Xin KPP random drift]}. $\square$

\

\textbf{Theorem 5.1.} (Large deviation principle for hitting time)
\textit{Almost surely with respect to $\Prob$ the following estimates hold.
Let $v, c\in \R$ and $c<v$. For any closed set $G \subset
(0,\infty)$ we have
$$\limsup\li_{t\ra\infty}\dfrac{1}{t}\ln \Prob^W
\left(\dfrac{T_{ct}^{vt}}{t}\in G\right)\leq -(v-c)\inf\li_{a\in
G}I^+\left(\dfrac{a}{v-c}\right) \
,$$$$\limsup\li_{t\ra\infty}\dfrac{1}{t}\ln \Prob^W
\left(\dfrac{T_{vt}^{ct}}{t}\in G\right)\leq -(v-c)\inf\li_{a\in
G}I^-\left(\dfrac{a}{v-c}\right) \ ;$$ and for any open set
$F\subset (0,\infty)$ we have
$$\liminf\li_{t\ra\infty}\dfrac{1}{t}\ln \Prob^W
\left(\dfrac{T_{ct}^{vt}}{t}\in F\right)\geq -(v-c)\inf\li_{a\in
G}I^+\left(\dfrac{a}{v-c}\right) \ , $$
$$\liminf\li_{t\ra\infty}\dfrac{1}{t}\ln \Prob^W
\left(\dfrac{T_{vt}^{ct}}{t}\in F\right)\geq -(v-c)\inf\li_{a\in
G}I^-\left(\dfrac{a}{v-c}\right) \ .$$}

\

\textbf{Proof.} The proof of this theorem is the same as the proof
of \cite[Theorem 2.3]{[Nolen-Xin KPP random drift]}. For the sake of
completeness we shall briefly repeat it here. We prove the first and
third bounds for example. The second and fourth estimates are the
same. Let $\lb\leq 0$. Let us consider the upper bound first. We
have, by Chebyshev's inequality, that

$$\begin{array}{l}
\play{\limsup\li_{t\ra\infty}\dfrac{1}{t}\ln\Prob^W\left(\dfrac{T_{ct}^{vt}}{t}<\al\right)}
\\
\play{\leq \limsup\li_{t\ra\infty}\dfrac{1}{t}\ln\Prob^W\left(e^{\lb
T_{c t}^{vt}}>e^{\lb \al t}\right)}
\\
\play{\leq - \lb \al + \limsup\li_{t\ra\infty}\dfrac{1}{t}\ln
q(ct,vt, \lb)}
\\
\play{=- \lb \al+(v-c)\mu^+(\lb) \ .}
\end{array}$$

Thus we see that
$$\limsup\li_{t\ra\infty}\dfrac{1}{t}
\ln\Prob^W\left(\dfrac{T_{ct}^{vt}}{t}<\al\right) \leq -\sup\li_{\lb
\leq 0}(\lb
\al-(v-c)\mu^+(\lb))=-(v-c)I^+\left(\dfrac{\al}{v-c}\right)
 \ ,$$ since $\lb_c^+=0$.

We now derive the lower bound. Let $u\in (0,\infty)$ and $\dt>0$.
Let $B_\dt(u)=(u-\dt, u+\dt)$ be the $\dt$-ball centered at $u$. Let
$\lb_u\leq 0$ be such that
$$I^+\left(\dfrac{u}{v-c}\right)=\sup\li_{\lb\leq 0}
\left(\lb\dfrac{u}{v-c}-\mu^+(\lb)\right)=\lb_u\dfrac{u}{v-c}-\mu^+(\lb_u)
\ .$$

Now we make use of a Cram\'{e}r's change of measure. Let

$$\dfrac{d\Prob^{W,u,t}}{d\Prob^W}=\dfrac{1}{S_{u,t}}e^{\lb_uT^{vt}_{ct}}\1_{T_{ct}^{vt}<\infty} \ ,$$

$$S_{u,t}=\E^W[e^{\lb_u T_{ct}^{vt}}\1_{T_{ct}^{vt}<\infty}] \ .$$

Then we get

$$\begin{array}{l}
\play{\Prob^W\left(\dfrac{T_{ct}^{vt}}{t}\in B_{\dt}(u)\right)}
\\
\play{\geq e^{-\lb_u ut-\dt t
|\lb_u|}\Prob^{W,u,t}\left(\dfrac{T_{ct}^{vt}}{t}\in
B_{\dt}(u)\right)\E^W[e^{\lb_u T_{ct}^{vt}}\1_{T_{ct}^{vt}<\infty}]}
\ .
\end{array}$$

One can show in the same way as in \cite[page 77]{[Zeitouni LDP
RWRE]} and \cite{[Taleb]}, that

$$\liminf\li_{t\ra\infty}\dfrac{1}{t}\ln\Prob^{W,u,t}\left(\dfrac{T_{ct}^{vt}}{t}\in B_{\dt}(u)\right)=0 \ . \eqno(5.2)$$

Suppose we already have (5.2), then we can conclude that we have

$$\begin{array}{l}
\play{\liminf\li_{t\ra\infty}\dfrac{1}{t}\ln\Prob^W\left(\dfrac{T_{ct}^{vt}}{t}\in
B_\dt(u)\right)}
\\
\play{\geq -\lb_u u-\dt |\lb_u|+(v-c)\mu^+(\lb_u)}
\\
\play{=(v-c)\left[\lb_u\dfrac{u}{v-c}-\mu^+(\lb_u)\right]-\dt|\lb_u|}
\\
\play{=(v-c)I^+\left(\dfrac{u}{v-c}\right)-\dt|\lb_u|}
\end{array}$$
which implies the lower bound.
 $\square$

\

\textbf{Theorem 5.2.} (Large deviation principle) \textit{Almost
surely with respect to $\Prob$ the following estimates hold. Let
$v\in \R$ and $\kp\in (0,1]$. For any closed set $G\subset
[0,\infty)$ we have}

$$\limsup\li_{t\ra\infty}\dfrac{1}{\kp t}\ln \Prob^W\left(\dfrac{vt-X^{vt}(\kp t)}{\kp t}
\in G\right)\leq -\inf\li_{c\in G}cI^+\left(\dfrac{1}{c}\right) \
,$$ \textit{and for any open set $F\subset[0,\infty)$ we have}

$$\liminf\li_{t\ra\infty}\dfrac{1}{\kp t}\ln \Prob^W\left(\dfrac{vt-X^{vt}(\kp t)}{\kp t}
\in F\right)\geq -\inf\li_{c\in F}cI^+\left(\dfrac{1}{c}\right) \
.$$

\textit{For any closed set $G\subset (-\infty, 0]$ we have }
$$\limsup\li_{t\ra\infty}\dfrac{1}{\kp t}\ln \Prob^W\left(\dfrac{vt-X^{vt}(\kp t)}{\kp t}
\in G\right)\leq -\inf\li_{c\in G}|c|I^-\left(\dfrac{1}{|c|}\right)
\ ,$$ \textit{and for any open set $F\subset(-\infty,0]$ we have}

$$\liminf\li_{t\ra\infty}\dfrac{1}{\kp t}\ln \Prob^W\left(\dfrac{vt-X^{vt}(\kp t)}{\kp t}
\in F\right)\geq -\inf\li_{c\in F}|c|I^-\left(\dfrac{1}{|c|}\right)
\ .$$

\

\textbf{Proof.} We show the first two estimates as an example. The
last two estimates are the same. We shall make use of the duality

$$\Prob^W\left(\dfrac{vt-X^{vt}(\kp t)}{\kp
t}>c\right)=\Prob^W\left(X^{vt}(\kp t)<vt-c\kp t\right)\leq
\Prob^W\left(\dfrac{\widehat{T}_{(v-c\kp)t}^{vt}}{t}<\kp\right) \
.$$

Here $\widehat{T}_r^s=\inf\{t\geq 0: X^s(t)\leq r\}$. We have
$$\widehat{T}_r^s=\1(Y_0\not\in I_0)\sm+T_{Y_\tau}^{Y_\sm}+\1(Y_{\widehat{T}_r^s}\not\in I_0)\widetilde{\tau} \ .$$

Here $\sm$ is the first time that the process $Y_t$, starting from
$Y_0\not \in I_0$, arrives at $I_0$; $\tau$ is the first time that
the process $Y_t$ arrives at the first branching point $K$ on $\Gm$
with $x$-coordinate $\geq r$; $\widetilde{\tau}$ is the first time
that the process $Y_t$, starting from $Y_\tau$, arrives at
$Y_{\widehat{T}_r^s}$. We note that by our Assumption 6 in
probability $1$ the distances $\rho(Y_{\sm}, (s,0))$,
$\rho(Y_{\tau}, (r,0))$ are finite. On the other hand, as
$t\ra\infty$, by Law of Large Numbers for stationary sequences we
have $\dfrac{\sm}{\widehat{T}_{rt}^{st}}\ra 0$ and
$\dfrac{\widetilde{\tau}}{\widehat{T}_{rt}^{st}}\ra 0$ almost
surely. Thus for fixed $r<s$ we have
$$\dfrac{T_{rt}^{st}}{\widehat{T}_{rt}^{st}}\ra 1 \text{ as } t\ra\infty \ . \eqno(5.3)$$
From here we have

$$\limsup\li_{t\ra\infty}\dfrac{1}{\kp
t}\ln\Prob^W\left(\dfrac{vt-X^{vt}(\kp t)}{\kp t}>c\right)\leq
\limsup\li_{t\ra\infty}\dfrac{1}{\kp t}\ln
\Prob^W\left(\dfrac{T_{(v-c\kp)t}^{vt}}{t}<\kp\right)\leq
-cI^+\left(\dfrac{1}{c}\right)$$ which proves the upper bound.

We now derive the lower bound. We have, for $0<\ve<1$,

$$\begin{array}{l}
\play{\Prob^W\left(\dfrac{vt-X^{vt}(\kp t)}{\kp t}\in
B_\dt(u)\right)}
\\
\play{=\Prob^W\left(X^{vt}(\kp t)\in B_{\kp t \dt}((v-\kp
u)t)\right)}
\\
\play{\geq \Prob^W\left(\widehat{T}_{(v-\kp u)t}^{vt}\in ((1-\ve)\kp
t, \kp t)\right)-\Prob^W\left(\sup\li_{(1-\ve)\kp t\leq s \leq \kp
t}|X^{vt}(s)-(v-\kp u)t|\geq \kp t \dt\right)}  \ .
\end{array}$$

The second term in the above formula can be estimated by using space
reversal invariance and the corresponding large deviation principle,
in the same way as \cite[proof of Theorem 2.4]{[Nolen-Xin KPP random
drift]} and \cite[Section 5]{[Taleb]}, provided that we have (5.3).
It turns out that

$$\lim\li_{\ve\ra 0}\limsup\li_{t\ra\infty}\dfrac{1}{t}\ln
\Prob^W\left(\sup\li_{(1-\ve)\kp t\leq s \leq \kp
t}|X^{vt}(s)-(v-\kp u)t|\geq \kp t \dt\right)=-\infty  \ .$$

So then we have, by (5.3) again and Theorem 5.1,

$$\begin{array}{l}
\play{\liminf\li_{t\ra\infty}\dfrac{1}{t}\ln\Prob^W\left(\dfrac{vt-X^{vt}(\kp
t)}{\kp t}\in B_\dt(u)\right)}
\\
\play{\geq \liminf\li_{\ve \ra
0}\liminf\li_{t\ra\infty}\dfrac{1}{t}\ln\Prob^W\left(T_{(v-\kp
u)t}^{vt}\in ((1-\ve)\kp t, \kp t)\right)}
\\
\play{=-\kp u I^+\left(\dfrac{1}{u}\right) \ .}
\end{array}$$
This proves the upper bound. $\square$

\section{Wave front propagation for reaction diffusion in narrow random channels}

After we get the quenched large deviation principle we study the
wave front propagation of the solution $u(t,(x,k))$ of (1.5) making
use of the arguments of \cite{[Nolen-Xin CMP]}, \cite{[Nolen-Xin KPP
random drift]} and \cite[Chapter 7]{[F red book]}.

We define non-random constants $c^*_+>0$ and $c^*_-<0$ as the
solutions of the equations
$$c^*_+I^+\left(\dfrac{1}{c^*_+}\right)=f'(0)  \ , \eqno(6.1.1)$$
$$|c^*_-|I^-\left(\dfrac{1}{|c^*_-|}\right)=f'(0)  \ . \eqno(6.1.2)$$

These solutions exist and are unique due to Lemma 5.3.

\

\textbf{Theorem 6.1.} \textit{For any closed set $F\subset (-\infty,
c^*_-)\cup (c_+^*, \infty)$ we have
$$\lim\li_{t\ra\infty}\sup\li_{c\in F}u(t, (ct, k))=0$$ almost
surely with respect to $\Prob$. For any compact set $K\subset(c^*_-,
c^*_+)$ we have
$$\lim\li_{t\ra\infty}\inf\li_{c\in K}u(t, (ct, k))=1$$ almost
surely with respect to $\Prob$.}

\

This theorem can be proved in the same way as \cite[Theorem 1.1,
Lemma 4.1, Lemma 4.2]{[Nolen-Xin KPP random drift]}. We shall
briefly sketch the proof here. We need a sequence of auxiliary
lemmas.

\

\textbf{Lemma 6.1.} \textit{For any closed set $F\subset (-\infty,
c_-^*)\cup(c_+^*, +\infty)$ we have
$$\lim\li_{t\ra\infty}\sup\li_{c\in F}u(t, (ct, k))=0$$ almost
surely.}

\

\textbf{Proof.} By the KPP condition and (1.5) we have

$$u(t, (ct, k))\leq \E^W_{(ct, k)}\left[\exp(f'(0)t)g(X^{ct}(t))\right] \ .$$

We notice that the support of the function $g(x)$ is a compact set
$U\subset (-\infty, \infty)$. Without loss of generality let us
assume that $U=B_{\dt}=(-\dt,\dt)$ for some $\dt>0$. Therefore we
have

$$\begin{array}{l}
u(t, (ct, k))
\\
\play{\leq \|g\|\exp(f'(0)t)\Prob^W_{(ct, k)}\left(-\dt \leq
X^{ct}(t)\leq \dt\right)}
\\
\play{=\|g\|\exp(f'(0)t)\Prob^W_{(ct,
k)}\left(c+\dfrac{\dt}{t}\geq\dfrac{ct-X^{ct}(t)}{t}\geq
c-\dfrac{\dt}{t}\right) \ .}
\end{array}$$

We apply Theorem 5.2 with $\kp=1$ and $v=c$. As $t\ra\infty$ we see
that, for $c>0$ such that
$f'(0)-cI^+\left(\dfrac{1}{c}\right)<-\ve<0$ and $\ve>0$; or for
$c<0$ such that $f'(0)-|c|I^-\left(\dfrac{1}{|c|}\right)<-\ve<0$ and
$\ve>0$ we have $\play{\limsup\li_{t\ra\infty}\dfrac{1}{t}\ln u(t,
(ct,k)) \leq -\dfrac{\ve}{2}}$. This proves the Lemma. $\square$

\

\textbf{Lemma 6.2.} \textit{For any compact set $K\subset (c^*_+,
+\infty)$ we have}

$$\liminf\li_{t\ra\infty}\dfrac{1}{t}\ln \inf\li_{c\in K}u(t, (ct,k))
\geq -\max\li_{c\in
K}\left[cI^+\left(\dfrac{1}{c}\right)-f'(0)\right] \ .
\eqno(6.2.1)$$

\textit{For any compact set $K\subset (-\infty, c^*_-)$ we have}

$$\liminf\li_{t\ra\infty}\dfrac{1}{t}\ln \inf\li_{c\in K}u(t, (ct,k))
\geq -\max\li_{c\in
K}\left[|c|I^-\left(\dfrac{1}{|c|}\right)-f'(0)\right] \ .
\eqno(6.2.2)$$

\

\textbf{Proof.} This lemma is proved in the same way as \cite[Lemma
4.1]{[Nolen-Xin KPP random drift]}, \cite[Corollary 1]{[Nolen-Xin
CMP]}, provided that we have the estimates (6.3) and (6.4) in the
following lemma. $\square$

\

\textbf{Lemma 6.3.} \textit{For any $v\in \R$ and $\eta>0$ we have}
$$\lim\li_{t\ra\infty}\sup\li_{|x|\leq |v|t}\Prob\left(\sup\li_{s\in [0,t]}|X^x(s)-x|\geq \eta t\right) =0 \ . \eqno(6.3)$$
\textit{Also, for a given $M>0$ there exists $\kp_0>0$ sufficiently
small so that} $$\limsup\li_{t\ra\infty}\sup\li_{|x|\leq
|v|t}\dfrac{1}{t}\ln \Prob\left(\sup\li_{s\in [0,\kp
t]}|X^x(s)-x|\geq \eta t\right)\leq -M \ , \eqno(6.4)$$\textit{
whenever $\kp<\kp_0$.}

\

\textbf{Proof.} This lemma is proved in the same way as \cite[Lemma
4.2]{[Nolen-Xin KPP random drift]}, by making use of Theorem 5.2. We
omit the details. $\square$

\

\textbf{Proof of Theorem 6.1.} With the above lemmas at hand the
lower bound follows from a standard argument as in \cite{[Nolen-Xin
CMP]} and \cite[Chapter 7, Theorem 3.1]{[F red book]}. We omit the
proof. $\square$

\

\end{document}